\documentclass[12pt]{article}
\usepackage{epsfig}

\title{ Fast Computation of the Economic Capital, the Value at Risk and the Greeks of a Loan Portfolio in the Gaussian Factor Model}

\author{ Pavel Okunev\footnote{This work was supported by the Director, Office of Science, Office  of Advanced Scientific Computing Research, of the U.S. Department of Energy under Contract No. DE-AC03-76SF00098.} \footnote{E-mail: pokunev@math.lbl.gov } \\  Department of Mathematics\\LBNL and UC Berkeley\\ Berkeley, CA 94720}

\date{July 1,2005}

\begin{document}
\maketitle

\begin{abstract}
We propose a fast algorithm for computing the economic capital, Value at Risk and Greeks in the Gaussian  factor model.   The algorithm proposed here is much faster than brute force Monte Carlo simulations or   Fourier transform based methods \cite{MD}. While the algorithm of Hull-White \cite{HW} is comparably fast, it assumes that all the loans in the portfolio have equal notionals and recovery rates. This is a very restrictive assumption which is unrealistic for many portfolios encountered in practice. Our algorithm  makes  no assumptions about the  homogeneity of the portfolio. Additionally, it is easier to implement than the algorithm of Hull-White. We use the implicit function theorem to  derive  analytic expressions  for the Greeks. 
\end{abstract}

\section{The Gaussian Factor Model}

Let us consider a portfolio of $N$ loans. Let the notional of loan $i$ be equal to $N_i$. Then the $i$th loan represents fraction
\begin{equation}
f_i=\frac{N_i}{\sum_j N_j}
\end{equation}
 of the notional of the whole portfolio. This means that if loan $i$ defaults  and the entire notional of the loan is lost the portfolio loses fraction $f_i$ or $100f_i\%$  of its value. In practice when a loan $i$ defaults a fraction $r_i$ of its notional will be recovered by the creditors.  Thus the actual loss given default (LGD) of loan $i$ is
\begin{equation}
LGD_i=f_i(1-r_i)
\end{equation}
fraction or
\begin{equation}
LGD_i=100f_i(1-r_i)\%
\end{equation}
of the notional of the entire portfolio.

We now describe the Gaussian m-factor model of portfolio losses from default. The model requires a number of input parameters. For each loan $i$ we are give a probability $p_i$ of its default. Also for each  $i$ and each $k=1,\ldots,m$ we are given a number $w_{i,k}$ such that $\sum_{k=1}^m w_{i,k}^2<1$. The number  $w_{i,k}$  is the loading factor of the loan $i$ with respect to factor $k$. Let $\phi_1, \ldots, \phi_m$ and $\phi^i, i=1,\ldots,N$ be independent standard normal random variables. Let $\Phi(x)$ be the cdf of the standard normal distribution.
In our model loan $i$ defaults if 
\begin{equation}
\sum_{k=1}^m w_{i,k}\phi_k+\sqrt{1-\sum_{k=1}^m w_{i,k}^2}\phi^i<\Phi^{-1}(p_i)
\end{equation}
This indeed happens with probability $p_i$.
  The factors $\phi_1,\ldots,\phi_m$ are usually interpreted  as the state of the global economy, the state of the  regional economy, the state of a particular industry and so on. Thus they are  the factors that affect the default behavior of all or at least a large group of loans in the portfolio. The factors $\phi^1,\ldots,\phi^N$ are interpreted as the idiosyncratic risks of the loans in the portfolio.
 
Let $I_i$ be defined by
\begin{equation}
I_i=I_{\{loan \ i\ defaulted\}}
\end{equation}
  We define the random loss caused by the default of loan $i$ as 
\begin{equation}
L_i=f_i(1-r_i)I_i,
\end{equation}
where $r_i$ is the recovery rate of loan $i$.
The total loss of the portfolio is 
\begin{equation}
L= \sum_i L_i
\end{equation}

An important property of the Gaussian factor model is that the $L_i$'s are not independent of each other. Their mutual  dependence is induced by the dependence of each $L_i$ on the common factors $\phi_1,\ldots,\phi_m$. Historical data supports the conclusion that losses due to defaults on different loans are correlated with each other. Historical data can also be used to calibrate the loadings $w_{i,k}$.the $L_i$'s are not independent of each other. Their mutual  dependence is induced by the dependence of each $L_i$ on the common factors $\phi_1,\ldots,\phi_m$. Historical data supports the conclusion that losses due to defaults on different loans are correlated with each other. Historical data can also be used to calibrate the loadings $w_{i,k}$.
 
\section{Conditional Portfolio Loss $L$}

When the values of the factors $\phi_1,\ldots,\phi_m$ are fixed, the probability of the default of loan $i$ becomes
\begin{equation}
p^i=\Phi\left( \frac{\Phi^{-1}(p_i)-\sum_kw_{i,k}\phi_k}{\sqrt{1-\sum_kw_{i,k}^2}} \right)
\end{equation}

The random losses $L_i$ become conditionally independent Bernoulli variables with the mean  given by 
\begin{equation}
E_{cond}(L_i)=f_i(1-r_i)p^i
\end{equation}
and the variance given by
\begin{equation}
VAR_{cond}(L_i)=f_i^2 (1-r_i)^2p^i(1-p^i)
\end{equation}   

By the Central Limit Theorem \cite{FE} the conditional distribution of the portfolio loss $L$, given the values of the factors $\phi_1,\ldots,\phi_m$,  can be approximated by the normal distribution with the mean
\begin{equation}
\label{m}
E_{cond}(L)=\sum_i E_{cond}(L_i)
\end{equation}
and the variance

\begin{equation}
\label{var}
VAR_{cond}(L)=\sum_i VAR_{cond}(L_i)
\end{equation}

We define the CDF of the unconditional portfolio loss $L$ by
\begin{equation}
F(x)=Prob(L\leq x), 
\end{equation}
for all real numbers $x$.
Since the conditional loss distribution is approximately normal, the CDF can be approximated by
\begin{equation}
\label{theformula}
F(x) \approx \int \Phi(E_{cond},VAR_{cond},x)\prod_{i=1}^{N} \rho(\phi_i)d\phi_1\ldots d\phi_m,
\end{equation}
where $\Phi(E_{cond},VAR_{cond},x)$ is the  CDF of the normal distribution with mean $E_{cond}$ and variance $VAR_{cond}$, while  $\rho(x)$ is the density of the standard normal distribution. 

\section{Numerical Example}
In this section we apply the proposed algorithm to the single factor Gaussian model of a portfolio with $N=125$ names. We choose a 125 name portfolio because it is the size of the standard DJCDX.NA.IG portfolio.  We choose a single factor model because it is the one most frequently used in practice.  The parameters of the portfolio are
\begin{eqnarray}
\label{port}
f_i&=&\frac{1}{N} \nonumber \\
p_i&=&0.015+\frac{0.05(i-1)}{N-1} \nonumber \\
r_i&=&0.5-\frac{0.1(i-1)}{N-1} \nonumber \\
w_{i1}&=&0.5-\frac{0.1(i-1)}{N-1},
\end{eqnarray}
where $i=1,\ldots,N$.

\begin{figure}[htbp]
\caption{CDF of the Unconditional Portfolio Loss}
\epsfig{file=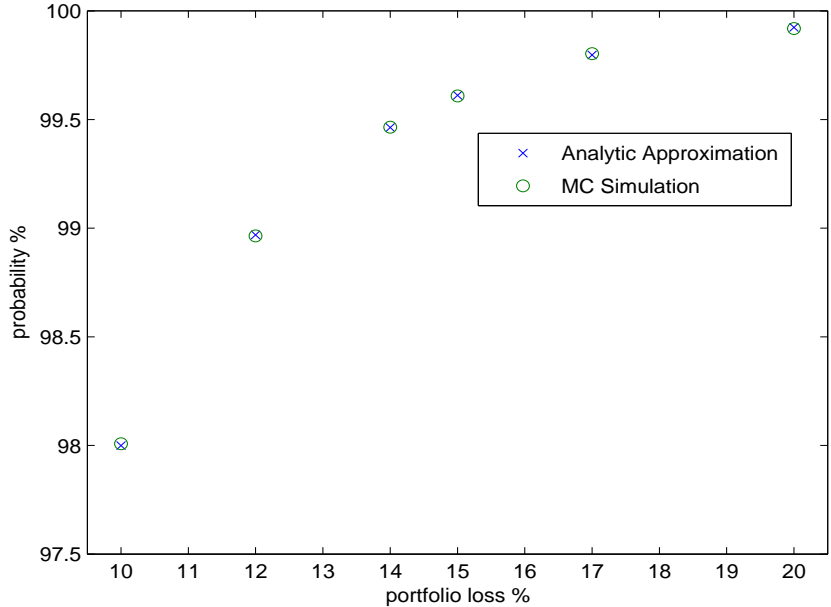,width=5.0 in,height=3.5 in}
\label{resultsfig}
\end{figure}

In Figure \ref{resultsfig} we compare the CDF  computed using $10^6$ Monte Carlo samples with the CDF computed using formula (\ref{theformula}).\footnote{The author has the code implementing the algorithm described here in MATLAB, VBA for Excel and C.} The agreement between the two is good. We add that the quality of approximation will improve  even further for a bigger portfolio and many bank portfolios have much more than 125 names.

\section{Computation of VaR and Economic Capital}
We define the $VaR$ (value at risk) of a given portfolio as the level of loss (expressed as fraction of portfolio notional), such that the probability of the portfolio loss being less or equal to $VaR$ is equal to a predefined confidence level $q$. Typically, $q$ is chosen to be between $0.99$ and $1$.
Thus we have 
\begin{equation}
Prob(L\leq VaR)=q
\end{equation}
From the results of the previous section  it follows that an accurate approximation to $VaR$  can be found by solving the equation
\begin{equation}
\label{equV}
\int \Phi(E_{cond},VAR_{cond},VaR)\prod_{i=1}^{N} \rho(\phi_i)d\phi_1\ldots d\phi_m=q
\end{equation}
This equation can be solved, for example, by bisection. To find the solution to (\ref{equV}) with accuracy of 1 basis point (0.01\%) we would need to evaluate the left hand side of (\ref{equV}) no more than 14 times.

{\bf Example.} {\em We want to calculate the $VaR$ of the portfolio (\ref{port}) with the confidence level $q=99.75\%$. We solve the equation (\ref{equV}) and round the solution to the nearest basis point to arrive at  $VaR=16.36\%$.  We now run Monte Carlo simulation with $5,000,000$ samples and compute the probability that the portfolio loss is less or equal to $16.36\%$. Rounded to the nearest basis point it turns out to be $99.75\%$. Thus the $VaR$ of the portfolio is indeed $16.36\%$.}

If desired, the convergence can be sped-up by using Newton's method.
After the $VaR$ is calculated we can calculate the economic capital by subtracting the average portfolio loss from $VaR$.

We now compare the algorithm proposed here with other proposed alternatives. The FFT based methods \cite{MD} require the computation of a large number  of Fourier transforms. To determine the $VaR$ with an error of less than 1 basis point (0.01\%) it is necessary to compute approximately 10,000 Fourier transforms.  Each Fourier transform is as expensive to evaluate as the left hand side of  (\ref{equV}). Thus our algorithm is significantly faster than the FFT based methods. It is well known that the FFT methods are  much faster than the direct Monte Carlo simulation. Thus our algorithm is much faster than the Monte Carlo approach. Finally, the recursive approach of Hull-White \cite{HW} is comparable in speed to the algorithm proposed here. However, it assumes that all the loans in the portfolio have equal notionals and recovery rates. This is a very restrictive assumption which is unrealistic for many portfolios encountered in practice. Our algorithm  makes  no assumptions about homogeneity of the portfolio. Additionally, it is easier to implement than the algorithm of Hull-White. 

\section{Computing the Greeks}
Since $VaR$ satisfies equation (\ref{equV}) we can use the implicit function theorem \cite{RD} to find its partial derivatives with respect to the parameters of the model. These partial derivatives are traditionally called the Greeks in finance.
We arrive at the following expressions
\begin{equation}
\frac{\partial}{\partial N_i}VaR=-\frac{\int \frac{\partial}{\partial N_i}\Phi(E_{cond},VAR_{cond},VaR)\prod_{i=1}^{N} \rho(\phi_i)d\phi_1\ldots d\phi_m}{\int \frac{\partial}{\partial VaR} \Phi(E_{cond},VAR_{cond},VaR)\prod_{i=1}^{N} \rho(\phi_i)d\phi_1\ldots d\phi_m}
\end{equation}

\begin{equation}
\frac{\partial}{\partial p_i}VaR=-\frac{\int \frac{\partial}{\partial p_i}\Phi(E_{cond},VAR_{cond},VaR)\prod_{i=1}^{N} \rho(\phi_i)d\phi_1\ldots d\phi_m}{\int \frac{\partial}{\partial VaR} \Phi(E_{cond},VAR_{cond},VaR)\prod_{i=1}^{N} \rho(\phi_i)d\phi_1\ldots d\phi_m}
\end{equation}

\begin{equation}
\frac{\partial}{\partial w_{i,k}}VaR=-\frac{\int \frac{\partial}{\partial w_{i,k}}\Phi(E_{cond},VAR_{cond},VaR)\prod_{i=1}^{N} \rho(\phi_i)d\phi_1\ldots d\phi_m}{\int \frac{\partial}{\partial VaR} \Phi(E_{cond},VAR_{cond},VaR)\prod_{i=1}^{N} \rho(\phi_i)d\phi_1\ldots d\phi_m}
\end{equation}

\begin{equation}
\frac{\partial}{\partial r_i}VaR=-\frac{\int \frac{\partial}{\partial r_i}\Phi(E_{cond},VAR_{cond},VaR)\prod_{i=1}^{N} \rho(\phi_i)d\phi_1\ldots d\phi_m}{\int \frac{\partial}{\partial VaR} \Phi(E_{cond},VAR_{cond},VaR)\prod_{i=1}^{N} \rho(\phi_i)d\phi_1\ldots d\phi_m}
\end{equation}
 and 
 \begin{equation}
\frac{\partial}{\partial q}VaR=\frac{1}{\int \frac{\partial}{\partial VaR} \Phi(E_{cond},VAR_{cond},VaR)\prod_{i=1}^{N} \rho(\phi_i)d\phi_1\ldots d\phi_m}
\end{equation}

 These expressions can be easily evaluated by numerical integration using Hermite-Gauss quadrature.

\section{Conclusions}
We proposed an algorithm for computing the $VaR$ and the economic capital of a loan portfolio in the Gaussian factor model.
The proposed method was tested on a portfolio of 125 names and gave high accuracy results. The accuracy will be even higher for portfolios with more names. Many of the bank portfolios are much larger than 125 names.
The proposed algorithm is much faster than the FFT based methods \cite{MD} and the brute force Monte Carlo approach. The speed of the Hull-White algorithm \cite{HW} is  comparable to that of the the algorithm proposed here, but  the Hull-White algorithm requires  that all the loans in the portfolio have equal notionals and recovery rates. This is a very restrictive assumption which is unrealistic for many portfolios encountered in practice. Our algorithm  makes  no assumptions about the homogeneity of the portfolio. Also, it is easier to implement than the algorithm of Hull-White.

Additionally, we obtained analytical expressions for the Greeks using the implicit function theorem. 

We also comment that the algorithm can be extended trivially to the case of stochastic  recovery rates and recovery rates correlated with the state of the factor variables.

Some of the ideas used in this paper were previously explored in \cite{PO} and \cite{PO2}.

\section{Acknowledgments}
I thank my adviser A. Chorin for his help and guidance during my time in UC Berkeley. I thank Mathilda Regan and Valerie Heatlie for their help in preparing this article. Last, but not least, I thank my family for their constant support.

\end{document}